\documentclass[reqno, a4paper]{amsart}

\usepackage{amsmath,amsfonts,amssymb}
\usepackage{verbatim}
\usepackage{enumerate}
\usepackage{url}
\usepackage{tikz}
\usepackage[left=2cm,right=2cm, top=3cm, bottom=2cm]{geometry}

\usetikzlibrary{arrows}

\def\N{\mathbb N}

\def\Z{\mathbb Z}

\def\ord{\mathop{\rm ord}\nolimits}

\theoremstyle{plain}
\newtheorem{theorem}{Theorem}[section]
\newtheorem{lemma}[theorem]{Lemma}
\newtheorem{definition}[theorem]{Definition}
\newtheorem{corollary}[theorem]{Corollary}
\newtheorem{proposition}[theorem]{Proposition}

\def\proof{\noindent {\it Proof: }}
\def\qed{\hfill\hbox{$\square$}}

\theoremstyle{definition}

\newtheorem*{Ack}{Acknowledgements}
\numberwithin{equation}{section}

\subjclass[2010]{20D60 (primary) and 11P70 (secondary)} 
\title{Extremal product-one free sequences in $C_q \rtimes_s C_m$}
\keywords{Zero-sum problem, Davenport constant, inverse zero-sum}

\author[F. E. Brochero Mart\'{\i}nez]{F. E. Brochero Mart\'{\i}nez}
\address{
Departamento de Matem\'{a}tica\\
Universidade Federal de Minas Gerais\\
UFMG\\
Belo Horizonte, MG\\
30123-970\\
Brazil\\
}
\email{fbrocher@mat.ufmg.br }

\author[S\'avio Ribas]{S\'avio Ribas}
\email{savio.ribas@gmail.com }

\date{\today}

\begin{document}

\maketitle

\begin{abstract}
Let $G$ be a finite group, written multiplicatively. The Davenport constant of $G$ is the smallest positive integer $d$ such that every  sequence of $G$ with $d$ elements has a non-empty subsequence with product $1$. Let $C_n \simeq \Z_n$ be the cyclic group of order $n$. In \cite{Bas}, J. Bass showed that the Davenport constant of the metacyclic group $C_q \rtimes_s C_m$, where $q$ is a prime number and $\ord_q(s) = m \ge 2$, is $m+q-1$. In this paper, we explicit the form of all sequences $S$ of $C_q \rtimes_s C_m$, with $q+m-2$ elements, that are free of product-$1$ subsequences.
\end{abstract}

\section{Introduction}

Given a finite group $G$ written multiplicatively, the {\em Zero-sum Problems} study conditions to ensure that a given sequence in $G$ has a non-empty subsequence with prescribed properties (such as  length, repetitions, weights) such that the {\em product} of its elements, in some order, is equal to the identity of the group. 

One of the first problem of this type is the remarkable Theorem of  Erd\"os-Ginzburg-Ziv (see \cite{EGZ}): Given $2n-1$ integers, it is possible to select $n$ of them, such that their sum is divisible by $n$, or in group theory language, every sequence $S$ with $l\ge 2n-1$ elements in a finite cyclic group of order $n$ has a subsequence of length $n$, the product of whose $n$ elements being the identity. 
In this theorem, the number $2n-1$ is the smallest integer with this property.

Traditionally, this class of problems have been extensively studied for abelian groups. We can see overviews of Zero-sum Theory for finite abelian groups in the surveys of Y. Caro \cite{Car} and W. Gao and A. Geroldinger \cite{GaGe}.

An important type of Zero-sum Problem is to determine the {\em Davenport constant} of a finite group $G$ (written multiplicatively): This constant, denoted by $D(G)$, is the smallest positive integer $d$ such that every sequence with $d$ elements in $G$ (repetition allowed) contains some subsequence such that the product of its terms in some order is $1$.

For $n \in \N$, let $C_n \simeq \Z_n$ denote the cyclic group of order $n$ written multiplicatively. The Davenport constant is known for some groups, such as:
\begin{itemize}
\item $D(C_n) = n$;
\item $D(C_m \times C_n) = m + n - 1$ if $m|n$  (J. Olson, \cite{Ols2});
\item $D(C_{p^{e_1}} \times \dots \times C_{p^{e_r}}) = 1 + \sum_{i = 1}^r (p^{e_i} - 1)$ (J. Olson, \cite{Ols1});
\item $D(D_{2n}) = n + 1$ where $D_{2n}$ is the Dihedral Group of order $2n$(see \cite{OlWh} and \cite{ZhGa});
\item $D(C_q \rtimes_s C_m) = m + q - 1$ where $q \ge 3$ is a prime number and $\ord_q(s) = m \ge 2$ (J. Bass, \cite{Bas}).
\end{itemize}
However, it is still open for most  other groups.

By the definition of the Davenport constant, there exist sequences $S$ of $G$ with $D(G) - 1$ elements that are free of {\em product-$1$ subsequences}, i.e, there exists $S=(x_1,\dots, x_{D(G)-1})$ sequence of $G$ such that $x_{i_1}x_{i_2}\cdots x_{i_k}\ne 1$ for every non empty subset $\{i_1,i_2,\dots,i_k\}\subset\{1,2,\dots,D(g)-1\}$ . The {\em Inverse Zero-sum Problems} study the structure of these extremal sequences which are free of product-$1$ subsequences with some prescribed property. Some overviews on the inverse problems can be found in  articles such as \cite{GaGeSc}, \cite{Sch} and \cite{GaGeGr}.

The inverse problems associated to the Davenport constant are already solved for a few abelian groups. The following theorem resolves the issue for the group $C_n$.

\begin{theorem}[{\cite[Theorem~2.1]{GaGeSc}}] \label{propC}
Let $S$ be a sequence in $C_n$ free of product-$1$ subsequences, where $n \ge 3$. Suppose that $|S| \ge (n+1)/2$. Then there exists some $g \in S$ with multiplicity $\ge 2|S| - n + 1$. In particular, if $|S| = n - 1 = D(C_n) - 1$ then $S = (\underbrace{g,\dots,g}_{n-1 \text{ times}})$, where $g$ is a generator of $C_n$.
\end{theorem}

Observe that in the cyclic case, sequences free of product-$1$ subsequences contain an element repeated many times. It is natural to ask if this is true in a general case. Specifically, let us say that a given finite abelian group $G$ has  {\em Property C} if every maximal sequence $S$ which is free of product-$1$ subsequences with at most $\exp(G)$ elements, i.e. the exponent of the group $G$,  has the form 
$$S = (\underbrace{T, T, \dots, T}_{\exp(G) - 1 \text{ times}}),$$
 for some subsequence $T$ of $S$. The above theorem states that $C_n$ has the Property C. It follows from a result of C. Reiher in \cite{Rei} that $C_p^2$ possesses Property C (see also \cite{GaGeSc} and \cite{GaGe3}). In \cite{GaGeGr}, W. Gao, A. Geroldinger and D. J. Grynkiewicz showed that this result is multiplicative, extending this result for $C_n^2$ where $n$ is a composite number. In \cite{Sch2}, W. A. Schmid discusses the case $C_n \times C_m$, where $n|m$. Not much is known about groups of rank $\ge 3$, only few specific cases (see, for example, \cite{Sch2}).

A {\em minimal zero sequence} $S$ in a finite abelian group $G$ is a sequence such that the product of its elements is $1$, but each proper subsequence is free of product-$1$ subsequences. In \cite[Theorem~6.4]{GaGe2}, W. Gao and A. Geroldinger showed that if $|S| = D(G)$ then $S$ contains some element $g \in G$ with order $\ord(g) = \exp(G)$ for certain groups such as $p$-groups, cyclic groups, groups with rank two and groups that are the sum of two elementary $p$-groups. They also conjectured that the same conclusion holds for every finite abelian group.

Another type of inverse zero-sum problem, associated to the Erd\"os-Ginzburg-Ziv Theorem, was proved by A. Bialostocki and P. Dierker in \cite{BiDi}. They established that if $S$ is a sequence in $C_n$ with $2n - 2$ elements and $S$ is free of product-$1$ subsequences with $n$ elements, then 
$$S = (\underbrace{g, \dots, g}_{n - 1 \text{ times}}, \underbrace{h, \dots, h}_{n - 1 \text{ times}})$$
 for some $g,h \in C_n$ with $\ord(gh^{-1}) = n$.

In this article, we caracterize the maximal sequences which are free of product-$1$ subsequences for certain non-abelian groups and we show that these sequences have a property similar to Property C. Let $q$ be a prime number, $m \ge 2$ be a divisor of $q-1$ and $s \in \Z_q^*$ such that $\ord_q(s) = m$. Denote by $C_q \rtimes_s C_m$ the {\em metacyclic group} $\Z_q \rtimes_s \Z_m$ written multiplicatively, i.e., the group generated by $x$ and $y$ with relations:
\begin{equation}\label{definicao}
x^m = 1, \quad y^q = 1, \quad yx = xy^s, \quad \text{ where } \ord_q(s) = m.
\end{equation}

Specifically, we prove the following result:

\begin{theorem}\label{main}
Let $q$ be a prime number, $m \ge 2$ be a divisor of $q-1$ and $s \in \Z_q^*$ such that $\ord_q(s) = m$, where $(m,q) \neq (2,3)$. Let $S$ be a sequence in the metacyclic group $C_q \rtimes_s C_m$ with $m+q-2$ elements. The following statements are equivalent:
\begin{enumerate}[(i)]
\item $S$ is free of product-$1$ subsequences;
\item For some $1 \le t \le q-1$, $1 \le i \le m-1$ such that $\gcd(i,m) = 1$ and $0 \le \nu_1, \dots, \nu_{m-1} \le q-1$,
$$S = ( \underbrace{y^t, y^t, \dots, y^t}_{q-1\text{ times}}, x^iy^{\nu_1}, x^iy^{\nu_2}, \dots, x^iy^{\nu_{m-1}} ).$$
\end{enumerate}
\end{theorem}

\section{Notation}

Let $G$ be a finite group written multiplicatively and $S = (g_1, g_2, \dots, g_l)$ be a sequence of elements of $G$. Suppose that $T$ is a subsequence of $S$, say $T = (g_{n_1}, g_{n_2}, \dots, g_{n_k})$, where $\{n_1, n_2, \dots, n_k\}$ is a subset of $\{1, 2, \dots, l\}$. We say that $T$ is a {\em product-$1$ subsequence} when
$$g_{\sigma(n_1)}g_{\sigma(n_2)}\dots g_{\sigma(n_k)} = 1$$
 for some permutation $\sigma$ of $\{n_1,\dots,n_k\}$, and if there are no such product-$1$ subsequences then we say that $S$ is {\em free of product-$1$ subsequences}.

Suppose that $S_1 = (g_{i_1},\dots, g_{i_u})$ and $S_2 = (g_{j_1}, \dots, g_{j_v})$ are subsequences of $S$. Then

\begin{itemize}
\item $\pi(S) = g_1 g_2 \dots g_l$ denotes the product of the elements in $S$ in the order that they appear;
\item $\pi_n(S) = g_{n+1} \dots g_l g_1 \dots g_n$, for $0 \le n \le l-1$, denotes the product of the elements in $S$ with a $n$-shift in the indices;
\item $|S| = l$ denotes the length of the sequence $S$;
\item $SS_1^{-1}$ denotes the subsequence formed by the elements of $S$ without the elements of $S_1$;
\item $S_1 \cap S_2$ denotes the intersection of the subsequences $S_1$ and $S_2$. In the case that $S_1 \cap S_2 = \varnothing$, we  say that $S_1$ and $S_2$ are disjoint subsequences;
\item $S_1S_2 = (g_{i_1}, \dots, g_{i_u}, g_{j_1}, \dots, g_{j_v})$ denotes the concatenation of $S_1$ and $S_2$;
\item $S^k = SS \dots S$ denotes the concatenation of $k$ identical copies of $S$'s.
\end{itemize}

For the group $C_q \rtimes_s C_m=\langle x,y| x^n=1,\  y^q=1,\  yx=xy^s\rangle$, let
\begin{itemize}
\item $H$ be the cyclic subgroup of order $q$ generated by $y$;
\item $N = (C_q \rtimes_s C_m) \setminus H = N_1 \cup N_2 \cup \dots \cup N_{m-1}$, where $N_i := x^iH$ for $1 \le i \le m-1$.
\end{itemize}
From the definition of semi-direct product,  $C_q \simeq H \vartriangleleft (C_q \rtimes_s C_m)$ and $C_m \simeq [(C_q \rtimes_s C_m) / H]$.

\section{Auxiliary results}

In this section we present the auxiliary theorems and lemmas that we use throughout the paper. First, we need the definition of sum-set and product-set:

\begin{definition}
If $X$ and $Y$ are non-empty subsets of an abelian group $G$ then the sum-set is defined by 
$$X + Y = \{a + b \in G \; | \;\; a \in X, b \in Y\}.$$

If $G$ is a non-abelian group, written multiplicatively, then the product-set is defined by 
$$X \cdot Y = \{a \cdot b \in G \; | \;\; a \in X, b \in Y\}.$$
\end{definition}

\vspace{0.3cm}

A very fundamental result on sum-sets is the Cauchy-Davenport Theorem, which gives a lower bound for the number of elements of a sum-set in $\Z_q$ depending on the cardinality of each set.

\begin{theorem}[Cauchy-Davenport inequality, {\cite[p.~44-45]{Nat}}]
For $q$ a prime number and for any $r$ non-empty sets $X_1, \dots, X_r \subset \Z_q$, 
$$|X_1 + \dots + X_r| \ge \min\{q, \; |X_1| + \dots + |X_r| - r + 1\}.$$
\end{theorem}

\vspace{0.3cm}

Looking at the inequality above in the case that $r = 2$ we get $|X + Y| \ge \min\{q, \; |X| + |Y| - 1\}$. A pair of subsets $X, Y$ of $\Z_q$ is called a {\em critical pair} if the equality $|X + Y| = \min\{q, \; |X| + |Y| - 1\}$ occurs. The following theorem provides  criteria for a pair $X, Y \subset \Z_q$ to be a critical pair.

\begin{theorem}[Vosper, \cite{Vos}]
Let $q$ be a prime number and let $X, Y$ non-empty subsets of $\Z_q$. Then 
$$|X + Y| = \min\{q, \; |X| + |Y| - 1\}$$
 if and only if one of the following conditions is satisfied:
\begin{enumerate}[(a)]
\item $|X| + |Y| > q$;
\item $\min\{|X|, |Y|\} = 1$;
\item $|X + Y| = q - 1$ and $Y = \Z_q \setminus \{c - a | \; a \in X\}$, where $\{c\} = \mathbb Z_q \setminus (X + Y)$;
\item $X$ and $Y$ are arithmetic progressions with the same common difference.
\end{enumerate}
\end{theorem}

Notice that  assertion $(a)$ above is the only one giving the equality $|X+Y| = q$, that is, $X+Y = \Z_q$. Assertion $(b)$ means that we just translate the set $Y$, supposing $|X| = 1$. The non-trivial cases yielding  a critical pair are $(c)$ and $(d)$.

\vspace{0.3cm}

Now, suppose that $A$ is a sequence in $C_q \rtimes_s C_m$ and that $A$ has no elements in the normal subgroup $H$, but that  the product of its elements is in $H$. 
The next lemma shows that if $A$ is  minimal with these properties then it is possible to generate at least $|A|$ distinct products in $H$, just shifting the order of the product.

\begin{lemma}\label{pinA}
Let $A = (a_1,a_2,\dots,a_l)$ be a sequence in $N \subset C_q \rtimes_s C_m$, where $\ord_q(s) = m$, such that $\pi(A) \in H$ but no subsequence of $A$ has product in $H$. Then:
\begin{enumerate}[(a)]
\item $\pi_n(A) \in H$ for every $0 \le n \le l-1$;
\item $\pi_i(A) \neq \pi_j(A)$ for every $0 \le i < j \le l-1$.
\end{enumerate}
\end{lemma}

\proof
This proof is contained in the proof of Lemma 14 in \cite{Bas}.
\qed

\vspace{0.3cm}

Throughout this paper, we deal with many expressions of the type $a_0 + a_1s + \dots + a_{m-1}s^{m-1} \pmod q$, where $a_0, a_1, \dots, a_{m-1} \in \Z_q$. The following lemma explains why we can assume without loss of generality that this kind of expression has a value different than $0$ under a weak assumption.

\begin{lemma}\label{wlog}
Let $A = (a_0, a_1, \dots, a_{m-1})$ be a sequence in $\Z_q$ with at least $2$ distinct elements, say $a_i \not\equiv a_j \pmod q$. Suppose that $\ord_q(s) = m$. Then 
$$\sum_{k=0}^{m-1} a_ks^k \not\equiv a_is^j + a_js^i + \sum_{k=0 \atop {k \neq i, j}}^{m-1} a_ks^k \pmod q.$$
 In particular, if one of them is $0$ modulo $q$ then the other is not $0$ modulo $q$.
\end{lemma}

\proof
Since at least two elements are distinct, there exists $0 \le i \le m-1$ such that $a_i \not\equiv a_{i+1} \pmod q$ (assuming that the index of the coefficients are taken modulo $m$, i.e., $a_0 = a_m$). Since $\ord_q(s) = m$, we may translate the coefficients by multiplying by $s$, therefore we may assume without loss of generality that $i = 0$, i.e., $a_0 \not\equiv a_1 \pmod q$. Now, if 
$$\sum_{k=0}^{m-1} a_ks^k \equiv a_0s + a_1 + \sum_{k=2}^{m-1} a_ks^k \pmod q$$
 then 
$$a_0 + a_1s \equiv a_1 + a_0s \pmod q \iff (a_0 - a_1)(1 - s) \equiv 0 \pmod q,$$
a  contradiction.
\qed

\vspace{0.3cm}

One of the assertions from Vosper's Theorem says that if $X$ and $Y$ are arithmetical progressions then $X$ and $Y$ form a critical pair. On the other hand, if $\alpha = a_0 + a_1s + \dots + a_{m-1}s^{m-1} \in \Z_q$ then the set $\{\alpha s^j\}_{j=0,1,\dots,m-1}$ is a geometric progression. Since we have the freedom to choose the order of the products, it is expected to deal with critical pairs that are formed by geometric progressions. The following lemma shows under some conditions that, in $\Z_q$, an arithmetic progression can not be a geometric progression.

\begin{lemma}\label{sinvariant}
Let $q \ge 5$ be a prime number, $s \in \Z_q^* \setminus \{1\}$ and $2 \le k \le q-1$. Let
\begin{align*}
\mathcal A &= \{1,2,3,\dots,k-1\}, \\
\mathcal B &= \{k,k+1,k+2,\dots,q-1\}
\end{align*}
be two sets of classes modulo $q$. Then $\mathcal A$ and $\mathcal B$ are not invariant by multiplication by $s$.
\end{lemma}

\proof
Since $\mathcal A$ and $\mathcal B$ are complementary in $\Z_q^*$, they are both $s$-invariants or not simultaneously. Suppose that they are $s$-invariants.

As $1\cdot s \in \mathcal A$ and $(q-1)\cdot s \in \mathcal B$, we obtain $2 \le s \le \min\{k-1,q-k\}$, therefore we may assume without loss of generality $k - 1 \le q -  k$. So $2 \le s \le k - 1 \le (q - 1)/2$. Let $c \equiv s^{\ord_q(s)-1} \pmod q$ be the inverse of $s$ modulo $q$. Since $s^j \in \mathcal A$ for all $j \in \N$ and $s \not\equiv 1 \pmod q$, we obtain $c \in \mathcal A$ and $c \not\equiv 1 \pmod q$, thus $c-1 \in \mathcal A$, which implies $s \cdot (c-1) \equiv 1 - s \in \mathcal A$. But $1 - s$ has a representative in $\mathcal B$, which is a contradiction.
\qed

\vspace{0.3cm}

The next lemma gives both upper and lower bounds for the number of solutions $(z,w) \in \Z_q^2$ of the equation $az^2 - bw^4 \equiv c \pmod q$ when $q \equiv 1 \pmod 4$ and $a,b \in \Z_q^*$ are fixed. We use the upper bound  to show that the set $\{c + bw^4 | \; w \in \Z_q^*\}$ contains both quadratic and non-quadratic residues modulo $q$.

\begin{lemma}\label{axbyc}
Let $q \equiv 1 \pmod 4$ be a prime number, $a, b, c \in \Z_q^*$ and $N$ the number of solutions $(z,w) \in \Z_q^2$ of $az^2-bw^4 \equiv c \pmod q$. Then
$$|N - q| < 3 \sqrt q.$$
\end{lemma}

\proof
Direct consequence of Theorem 5 page 103 in \cite{IrRo}.
\qed

\vspace{0.3cm}

\begin{corollary}\label{squarenonsquare}
Let $q \equiv 1 \pmod 4$ be a prime number such that $q \ge 13$ and let $b,c \in \Z_q^*$ . Then there exist $w_1, w_2\in \Z_q^*$ such that $c + bw_1^4$ is a quadratic residue and $c + bw_2^4$ is a non-quadratic residue. 
\end{corollary}

\proof
Fix $a \in \Z_q^*$ and denote by $N$ the number of solutions of 
\begin{equation}\label{equacaoabc}
az^2 - bw^4 \equiv c \pmod q.
\end{equation}

Let $w_0 \in \Z_q^*$ such that $\ord_q(w_0) = 4$. Since each solution of equation (\ref{equacaoabc}) generates seven other solutions (by switching $z$ by $-z$ and multiplying $w$ by powers of $w_0$), the number of elements in the set $$\mathcal C=\{az^2 | \; z \in \Z_q^*\}\cap \{ c + bw^4 | \; w \in \Z_q^*\}$$ is at most $N/8$. By the previous lemma, it follows that 
$$|\mathcal C| < \frac {q+3\sqrt q}8 < \frac{q-1}4 = |\{c + bw^4 | \; w \in \Z_q^*\}|,$$
 therefore 
$$\{c + bw^4 | \; w \in \Z_q^*\} \not\subset \{az^2 | \; z \in \Z_q^*\}.$$

Thus, selecting $a$ being either a square or not, we obtain that the set $\{ c + bw^4 | \; w \in \Z_q^*\}$ contains quadratic residues and non-quadratic residues.
\qed

\vspace{0.3cm}

It follows from the corollary above and Lemma \ref{wlog} that, in the case $m = (q-1)/2$, it is possible to get more than the $m$ distinct values which were provided by  Lemma \ref{pinA}.

\begin{corollary}\label{q=1,r=2}
Let $q \equiv 1 \pmod 4$ be a prime number such that $q \ge 13$, $s \in \Z_q^*$, such that $\ord_q(s) = \frac{q-1}2 = m$ and $a_0, a_1,\dots a_{m-1} \in \Z_q$. Suppose that each of the sets 
$$\{a_0,a_2,a_4,\dots,a_{m-2}\} \text{ and } \{a_1,a_3,a_5,\dots,a_{m-1}\}$$
 has at least two distinct elements modulo $q$. Then the set 
$$\mathcal A = \left\{a_{\sigma(0)} + a_{\sigma(1)}s + a_{\sigma(2)}s^2 + \cdots + a_{\sigma(m-1)}s^{m-1} \in \Z_q \; | \; \sigma \text{ is a permutation of } (0, 1, \dots, m-1) \right\}$$
  has at least $q-1$ distinct elements. 
\end{corollary}

\proof
By Lemma \ref{wlog}, we can suppose without loss of generality that 
$$\alpha := a_0 + a_1s + a_2s^2 + \cdots + a_{m-1}s^{m-1} \not\equiv 0 \pmod q.$$

Notice that we can obtain $\alpha s^j$ by shifting the coefficients, so $\alpha s^j \in \mathcal A$ for all $0 \le j \le m-1$. By Lemma \ref{pinA}, these $m$  elements are distinct.

Since $\ord_q(s) = (q-1)/2$, $s$ is a square modulo $q$ and, indeed, $s$ generates the quadratic residues modulo $q$, therefore the elements $\alpha s^j$ are all quadratic residues or all non-quadratic residues depending on whether $\alpha$ is a square or not. Define
\begin{align*}
b &\equiv a_1s + a_3s^3 + \cdots + a_{m-1}s^{m-1} \pmod q, \\
c &\equiv a_0 + a_2s^2 + \dots + a_{m-2} s^{m-2} \pmod q.
\end{align*}
In the same way, by Lemma \ref{wlog} we can assume without loss of generality that $b,c \not\equiv 0 \pmod q$. Since $c + bs^{2j} \in \mathcal A$ for all $0 \le j \le m/2$,  Corollary \ref{squarenonsquare} tells us that the set $\{c + bs^{2j} | \; 0 \le j \le m/2\}$ contains a quadratic residue and a non-quadratic residue modulo $q$. By multiplying by $s$, we obtain all quadratic residues and all non-quadratic residues modulo $q$. Therefore, $|\mathcal A| \ge q-1$.
\qed

\vspace{0,3cm}

Analogously to the previous corollary but now in the case $m = q-1$, the next result states that it is possible to split the coefficients into  parts such that each part generates at least $m$ distinct values, also improving the result of Lemma \ref{pinA}.

\begin{corollary}\label{q=1,r=1}
Let $q \equiv 1 \pmod 4$ be a prime number such that $q \ge 13$, $s$ be a generator of $\Z_q^*$ and $a_0, a_1,\dots a_{q-2} \in \Z_q$. Suppose that each of the sets 
$$\{a_0, a_4, a_8, \dots, a_{q-5}\}, \{a_1, a_5, a_9, \dots, a_{q-4}\}, \{a_2, a_6, a_{10}, \dots, a_{q-3}\} \text{ and } \{a_3, a_7, a_{11}, \dots, a_{q-2}\}$$
 has at least two distinct elements modulo $q$. Then each of the sets
\begin{align*}
\mathcal A_e &= \left\{a_{\sigma(0)} + a_{\sigma(2)}s^2 + a_{\sigma(4)}s^4 + \cdots + a_{\sigma(q-3)}s^{q-3} \in \Z_q \; | \; \sigma \text{ is a permutation of } (0, 2, 4, \dots, q-3) \right\}, \\
\mathcal A_o &= \left\{a_{\tau(1)}s + a_{\tau(3)}s^3 + a_{\tau(5)}s^5 + \cdots + a_{\tau(q-2)}s^{q-2} \in \Z_q \; | \; \tau\text{ is a permutation of } (1, 3, 5, \dots, q-2) \right\}
\end{align*}
have at least $q-1$ distinct elements. 
\end{corollary}

\proof
This  follows directly from the previous corollary.
\qed

\section{Sequences in $C_5 \rtimes_s C_m$}

In our proof of Theorem \ref{main}, for the case $q \equiv 1 \pmod 4$, we use Corollaries \ref{q=1,r=2} and \ref{q=1,r=1}, where the hypothesis $q \ge 13$ is necessary. In this section, we consider the remaining case, i.e., $q = 5$. There are two possibilities for $m$, namely, $m=2$ and $m=4$.

For $m=2$ the only possible value for $s$ is $4 \pmod 5$, therefore we obtain the Dihedral Group of order $10$. The next proposition deals with this case and shows that if a sequence $S$ is free of product-$1$ subsequences and satisfies some assumptions then $S$ must contain an element in $H$, the normal subgroup.

\begin{proposition}\label{mq25}
Let $S$ be a sequence in $C_5 \rtimes_4 C_2$ free of product-$1$ subsequences and suppose that $|S| = 5$. Then $S \cap H \neq \varnothing$.
\end{proposition}

\proof
Suppose that $S = (xy^{\alpha_1},\dots,xy^{\alpha_5})$, where $0 \le \alpha_1 \le \dots \le \alpha_5 \le 4$. If there exist two identical elements in $S$ then their product is $1$ and so $S$ is not free of product-$1$ subsequences. Thus, $(\alpha_1,\dots,\alpha_5) = (0,1,2,3,4)$, therefore $x \cdot xy^4 \cdot xy \cdot xy^2 = 1$. Hence, $S$ is not free of product-$1$ subsequences, a  contradiction.
\qed

\vspace{0,3cm}

For $m=4$, the cases to consider are $s \equiv 2 \pmod 5$ or $s \equiv 3 \pmod 5$. The proposition below shows that if $S$ is free of product-$1$ subsequences then $S$ can not belong to a single coset $N_i$, where $\gcd(i,4) = 1$.

\begin{proposition}\label{mq45}
Let $s \in \{2,3\}$ and let $S$ be a sequence in $C_5 \rtimes_s C_4$ such that $|S| = 7$. If every element of $S$ belongs to a coset $N_i$, where $i \in \{1,3\}$, then $S$ is not free of product-$1$ subsequences.
\end{proposition}

\proof
Let $S = (x^iy^{\alpha_1},\dots,x^iy^{\alpha_7})$, where $i \in \{1,3\}$ and $0 \le \alpha_1 \le \dots \le \alpha_7 \le 4$. Notice that at most three of the $\alpha_j$'s are equal, otherwise the product of four identical elements would be $1$. This implies that there are at least three distinct elements. Also, there are at least a pair of elements repeating two or three times. Since $\ord_5(s) = 4$, it follows that $s^{3i_0} + s^{2i_0} + s^{i_0} + 1 \equiv 0 \pmod 5$. Therefore, for all $\beta \in \Z_q$ it holds that:
\begin{align*}
x^{i_0}y^{a_1} \cdot x^{i_0}y^{a_2} \cdot x^{i_0}y^{a_3} \cdot x^{i_0}y^{a_4} &= y^{a_1s^{3i_0} + a_2s^{2i_0} + a_3s^{i_0} + a_4} \\
&= y^{a_1s^{3i_0} + a_2s^{2i_0} + a_3s^{i_0} + a_4 - \beta(s^{3i_0} + s^{2i_0} + s^{i_0} + 1)} \\
&= x^{i_0}y^{a_1 - \beta} \cdot x^{i_0}y^{a_2 - \beta} \cdot x^{i_0}y^{a_3 - \beta} \cdot x^{i_0}y^{a_4 - \beta},
\end{align*}
thus we may assume without loss of generality that the element that repeats most  is $\alpha_1 = 0$. If there is another pair of identical elements, say $1 \le \alpha_j = \alpha_{j+1} = \lambda \le 4$ for some $3 \le j \le 6$, then we may choose $A_1 = (x, xy^{\lambda}, x, xy^{\lambda})$ or $A_1 = (x^3, x^3y^{\lambda}, x^3, x^3y^{\lambda})$. Since $\ord_5(s) = 4$ and $\gcd(i,4) = 1$, 
it follows that $\lambda s^{2i} + \lambda \equiv 0 \pmod 5$, and so $x^i \cdot x^iy^{\lambda} \cdot x^i \cdot x^iy^{\lambda} = 1$. Otherwise, the only possibility is $(\alpha_1,\dots,\alpha_7) = (0,0,0,1,2,3,4)$ and we may choose, for example, 
\begin{align*}
A_1 &= (x,x,xy,xy^2) \text{ or } A_1 = (x,x,xy,xy^3) \text{ when } i = 1, \\
A_1 &= (x^3, x^3, x^3y, x^3y^2) \text{ or } A_1 = (x^3, x^3, x^3y, x^3y^3) \text{ when } i = 3,
\end{align*}
as the following table shows:

\begin{center}
\begin{tabular}{c||c|c}
 & $s=2$ & $s=3$ \\ \hline \hline
$i=1$ & $s^i \equiv_5 2 \Rightarrow x \cdot x \cdot xy \cdot xy^3 = 1$         & $s^i \equiv_5 3 \Rightarrow x \cdot x \cdot xy \cdot xy^2 = 1$ \\ \hline
$i=3$ & $s^i \equiv_5 3 \Rightarrow x^3 \cdot x^3 \cdot x^3y \cdot x^3y^2 = 1$ & $s^i\equiv_5 2 \Rightarrow x^3 \cdot x^3 \cdot x^3y \cdot x^3y^3 = 1$
\end{tabular}
\end{center}

\vspace{0,3cm}

Thus, $S$ is not free of product-$1$ subsequences.

\section{Proof of Theorem \ref{main}}\label{proofmain}

In this section we investigate the sequences of $C_q \rtimes_s C_m$ with $m+q-2$ elements which are free of product-$1$ subsequences and we prove that an analogue of Property C holds for $C_q \rtimes_s C_m$ provided $\ord_q(s) = m$ and $(m,q) \neq (2,3)$. Notice that if $(m,q) = (2,3)$ then $s \equiv 2 \pmod 3$ and the assertion of Theorem \ref{main} is not true for the group $C_3 \rtimes_2 C_2$, which is isomorphic to $D_6$, the Dihedral Group of order $6$, and to $S_3$, the Permutation Group of $3$ elements. In fact, the sequence $S = (x,xy,xy^2)$ provides the only counter example when $(m,q) = (2,3)$.

It is easy to check that $(ii)$ implies $(i)$, therefore we just need to prove that $(i)$ implies $(ii)$. Let us assume that $S$ is a sequence free of product-$1$ subsequences in $C_q \rtimes_s C_m$. Let us also define $k \in \Z$ by the equation $$|S \cap H| = q - k.$$

If $k \le 0$, since $D(H) = D(C_q) = q$, then there exists a product-$1$ subsequence in $H$.

If $k = 1$ then $|S \cap H| = q - 1$ and $|S \cap N| = m-1$. By Theorem \ref{propC}, the elements of $S \cap H$ must all be equal, say, $S \cap H = \{y^t\}^{q-1}$ and the other elements of $S$ must be in the same $N_i$, where $\gcd(i,m)=1$.

From now on, assume $k \ge 2$. In this case, we are going to prove that $S$ is not free of product-$1$ subsequences. We have 
$$|S \cap N| = m + k - 2 \ge m = D(C_m) = D((C_q \rtimes_s C_m)/H).$$

Let $A = (a_1,\dots,a_l)$ be a subsequence of $S \cap N$ such that $\pi(A) \in H \simeq C_q$ but no subsequence of $A$ has product in $H$, as in Lemma \ref{pinA}. We have that $2 \le |A| \le m$. Let us denote $A$ by $A_1$. If $|(S \cap N)A_1^{-1}| \ge m$ then we can construct $A_2$ with the same property and repeat this argument replacing successively $(S \cap N)(A_1\dots A_j)^{-1}$ by $(S\cap N)(A_1\dots A_jA_{j+1})^{-1}$ and so on, until $|(S\cap N)(A_1\dots A_r)^{-1}| \le m-1$. This implies that
$$\sum_{i=1}^r |A_i| \ge |S \cap N| - (m-1) = k-1.$$

We construct the following set of products:
\begin{align*}
R = (\{\pi_j(A_1)\}_{j=0,1,\dots, |A_1|-1}) \cdot &(\{1\} \cup \{\pi_j(A_2)\}_{j=0,1,\dots, |A_2|-1})\cdot  \dots \\
&\dots \cdot (\{1\} \cup \{\pi_j(A_r)\}_{j=0,1,\dots,|A_r|-1}) \subset H.
\end{align*}

By the Cauchy-Davenport inequality we obtain 
$$|R| \ge \min\left\{q, \; |A_1| + \sum_{i=2}^r (|A_i|+1) -  r + 1\right\} = \min\left\{q, \;\; \sum_{i=1}^r |A_i| \right\}.$$

If this minimum is $q$ then $R = H \ni 1$ and $S$ is not free of product-$1$ subsequences. On the other hand, if this minimum is $\sum_{i=1}^r |A_i| \ge k-1$ then we obtain at least $k-1$ distinct elements in $H$ arising from $S \cap N$. Suppose without loss of generality that $S \cap H = (h_1,h_2,\dots,h_{q-k})$. If $\sum_{i=1}^r |A_i| \ge k$ then, by the Pigeonhole Principle, either $1 \in R$ or $R$ contains the inverse of one of the products $h_1, (h_1h_2), \dots, (h_1\dots h_{q-k})$, hence $S$ is not free of product-$1$ subsequences.

Therefore, suppose that $\sum_{i=1}^r |A_i| = k-1$, namely, $R = \{ g_1, g_2, \dots, g_{k-1}\}$. If there exist $1 \le n_1 \le k-1$ and $1 \le n_2 \le q-k$ such that $g_{n_1} = (h_1\dots h_{n_2})^{-1}$ then $g_{n_1}h_1\dots h_{n_2} = 1$ and so $S$ is not free of product-$1$ subsequences. Therefore, 
$$\{g_1, \dots, g_{k-1}, h_1^{-1}, (h_1h_2)^{-1}, \dots, (h_1\dots h_{q-k})^{-1}\} = \{y, y^2, \dots, y^{q-1}\}.$$

If $h_i \neq h_j$ for some $1 \le i < j \le q-k$, say without loss of generality that $h_1 \neq h_2$, then either the set 
$$\{g_1, \dots, g_{k-1}, h_1^{-1}, h_2^{-1}, (h_1h_2)^{-1}, \dots, (h_1\dots h_{q-k})^{-1}\}$$
 has $q$ elements, in particular, it has the element $1$, or it contains two identical elements. In any case, $S$ has a product-$1$ subsequence. Hence, 
 $$h_1 = h_2 = \dots = h_{q-k} = y^t$$
  for some $1 \le t \le q-1$ and 
$$R = \{g_1, g_2, \dots, g_{k-1}\} = \{y^t, y^{2t}, \dots, y^{(k-1)t}\}.$$

Let $C := (S\cap N) (A_1\dots A_r)^{-1}$. Notice that $C$ is free of subsequences with product in $H$. Since 
$$|C| = |S \cap N| - \sum_{j=1}^r |A_j| = m+k-2 - (k-1) = m-1,$$
 we conclude, by Theorem \ref{propC}, that $C$ does not have subsequences with product in $H$ if these $m-1$ elements are in the same class $N_{i_0}, 1 \le i_0 \le m-1, \gcd(i_0,m)=1$. Thus, we assume that $C$ is a sequence in $S \cap N_{i_0}$, namely, 
 $$C = (x^{i_0}y^{t_1}, x^{i_0}y^{t_2}, \dots, x^{i_0}y^{t_{m-1}}).$$

If there exist $1 \le j \le r$ such that $A_j$ has at least one element out of $N_{i_0}$, then we may select $x^{i_1}y^{\theta} \in A_j$ with $i_1 \neq i_0$. If $2 \le l \le m-1$ is such that $i_0 l \equiv i_1 \pmod m$ then we may replace $A_j$ by 
$$\tilde{A_j} = (A_j \setminus \{x^{i_1}y^{\theta}\}) \cup (x^{i_0}y^{t_1}, \dots, x^{i_0}y^{t_l}),$$
 which also has product in $H$, since 
$$\prod_{n = 1}^l x^{i_0}y^{t_n} = x^{i_1}y^T$$ 
for some $T \in \Z_q$ and $H \vartriangleleft C_q \rtimes_s C_m$. Since $|\tilde{A_j}| > |A_j|$, the new set $R$ generated by this change has more elements. This replacement may mean that  $\tilde{A_j}$ is not  minimal anymore and in this case we break $\tilde{A_j}$ into its minimal components. Thus, 
 $$\sum_{n=1 \atop n \neq j}^r |A_n| + |\tilde{A_j}| > \sum_{n=1}^r |A_j| = k - 1,$$
  therefore there exists a product-$1$ subsequence in $S$.

Hence, $S \cap N$ must be a sequence in $N_{i_0}$, so 
$$A_j \!\!\!\! \pmod H = \{x^{i_0}\}^m \;\;\;\; \text{ for all }1 \le j \le r, \text{ and } \;\;\;\; C = (x^{i_0}y^{t_1}, \dots, x^{i_0}y^{t_{m-1}}).$$

By double counting the number of elements in $S \cap N_{i_0}$ we conclude that $m + k - 2 \equiv m - 1 \pmod m$, that is, $k \equiv 1 \pmod m$. Since $k \ge 2$ and $k \equiv 1 \pmod m$, we have $m + 1 \le k \le q$. 

Observe that the case $m+1 \le k < q$ is not possible. In fact, in this case, if 
$$A_j = (x^{i_0}y^{\eta_1}, \dots, x^{i_0}y^{\eta_m})$$
 then 
$$\pi(A_j) = y^{\eta_1 s^{i_0(m-1)} + \eta_2 s^{i_0(m-2)} + \dots + \eta_{m-1} s^{i_0} + \eta_m}$$
 and, more generally, 
$$\pi_n(A_j) = \left(y^{\eta_1 s^{i_0(m-1)} + \eta_2 s^{i_0(m-2)} + \dots + \eta_{m-1} s^{i_0} + \eta_m}\right)^{s^{i_0n}}.$$

We claim that $R$ is invariant under taking powers of $s^{i_0}$. In fact, since $\gcd(i_0,m) = 1$ we obtain $\ord_q(s^{i_0}) = \ord_q(s) = m$. An element in $R$ is of the form 
$$y^{jt} = \pi_{j_1}(A_{\nu_1}) \dots \pi_{j_u}(A_{\nu_u}),$$
 where $1 \le j \le k-1$. Taking powers of $s^{i_0}$ in both sides, we obtain that 
$$y^{s^{i_0}jt} = \pi_{j_1+1}(A_{\nu_1}) \dots \pi_{j_u+1}(A_{\nu_u})$$
 belongs to $R$, because it can be obtained by the product of some $A_i$'s in some order. Therefore, the claim is proved.

Looking at the exponent of $y$, the above claim implies that the set $\{t,2t,\dots,(k-1)t\}$ is $s^{i_0}$-invariant modulo $q$, and so $\{1,2,\dots,k-1\}$ is $s^{i_0}$-invariant, which contradicts  Lemma \ref{sinvariant}.

From now on, we assume $k = q$ 
and $S$ is a sequence in $N_{i_0}$. As $|A_j| = m$ for $1 \le j \le r$ and $|C| = m-1$ we obtain $mr + m - 1 = |S| = m + q - 2$, therefore $mr = q-1$. We consider the following cases depending on whether $r \ge 3$, $r = 2$ or $r = 1$:

\begin{enumerate}[(1)]
\item {\bf Case $r \ge 3$:} This implies that $3m \le q-1$. Since the equality in the Cauchy-Davenport inequality occurs, Vosper's Theorem with the sets 
$$\tilde{A} = \{\pi_n(A_1)\}_n \;\;\; \text{ and } \;\;\; \tilde{B} = (\{1\} \cup \{\pi_n(A_2)\}_n) \cdot \dots \cdot (\{1\} \cup \{\pi_n(A_r)\}_n)$$
 says that at least one of the following statements hold:
\begin{enumerate}[(i)]
\item $|\tilde{A}| + |\tilde{B}| > q$;
\item $\min\{|\tilde{A}|, |\tilde{B}|\} = 1$;
\item $|\tilde{A} \cdot \tilde{B}| = q - 1$ and $\tilde{A} = H \setminus \{b^{-1} | \; b \in \tilde{B}\}$;
\item $\tilde{A}$ and $\tilde{B}$ are arithmetic progressions with the same common difference.
\end{enumerate}

Since $|\tilde{A}| + |\tilde{B}| = q$, (i) is not possible, and since $\min\{|\tilde{A}|, |\tilde{B}|\} \ge m \ge 2$, (ii) does not hold. In order to discard item (iii), define
\begin{align*}
\tilde{B_1} &= \{1\} \cup \{\pi_n(A_2)\}_n, \\
\tilde{B_2} &= (\{1\} \cup \{\pi_n(A_3)\}_n) \cdot \dots \cdot (\{1\} \cup \{\pi_n(A_r)\}_n).
\end{align*}
As $3m \le mr = q-1$, we have
\begin{align*}
q - m &= |\tilde{B}| = |\tilde{B_1} + \tilde{B_2}| = |\tilde{B_1}| + |\tilde{B_2}| - 1 < q - 1, \\
|\tilde{B_1}| &= m+1, \\
|\tilde{B_2}| &= q - 2m \ge m+1.
\end{align*}

Therefore, the items (i), (ii) and (ii) from Vosper's Theorem are false, hence the exponents of the elements in each of the sets $\tilde{B_1}$ and $\tilde{B_2}$ form arithmetic progressions, say
$$\tilde{B_1} = \{y^{-av}, \dots, y^{-v}, 1, y^v, \dots, y^{(m-a-1)v}\}.$$

Looking at the Vosper's equality involving $\tilde{A}$ and $\tilde{B}$, item (iii) tells us that $\tilde{A}$ also form an arithmetic progression in the exponent with the same common difference, say
$$\tilde{A} = \{y^{w}, y^{w+v}, \dots, y^{w + (m-1)v}\}.$$
By switching $A_1$ and $A_2$, the only possibility is that $\{\pi_n(A_1)\} = \{y^v, y^{2v}, \dots, y^{mv}\}$ for some $1 \le v \le q-1$.

On the other hand, the exponents of the elements from the set $\{\pi_n(A_1)\}_n$ are $s^{i_0}$-invariant, that is, the exponents of $y$ in $R$ are invariant by multiplication by $s^{i_0}$. By Lemma \ref{sinvariant}, the set $\{\pi_n(A_1)\}_n$ can not be of the above form.

\item {\bf Case $r = 2$:} In this case, $m = (q-1)/2$ and, in particular, $s^{i_0}$ generates the quadratic residues modulo $q$. The case $(m,q) = (2,5)$ follows from Proposition \ref{mq25}, therefore we may assume $q \ge 7$. If there exist $m$ identical elements then their product is $1$, thus there are at most $m-1$ identical elements. Since 
$$\frac{m+q-2}{m-1} = \frac{3m-1}{m-1} > 3,$$
 there are at least $4$ distinct elements among $S = \{x^{i_0}y^{\alpha_1}, \dots, x^{i_0}y^{\alpha_{3m-1}}\}$. We split this case into the two subcases $q \equiv 1 \pmod 4$ and $q \equiv 3 \pmod 4$.

\begin{enumerate}[({2.}1)]
\item {\bf Subcase $q \equiv 1 \pmod 4$ and $q \ge 13$:} We may choose
\begin{align*}
A_1 &= (x^{i_0}y^{a_0}, x^{i_0}y^{a_1}, \dots, x^{i_0}y^{a_{m-1}}), \\
A_2 &= (x^{i_0}y^{b_0}, x^{i_0}y^{b_1}, \dots, x^{i_0}y^{b_{m-1}})
\end{align*}
to be disjoint subsequences of $S$ such that it is possible to split each one into two subsequences of the same size, say
\begin{align*}
A_1 &= (x^{i_0}y^{a_0}, x^{i_0}y^{a_2}, \dots, x^{i_0}y^{a_{m-2}}) \cup (x^{i_0}y^{a_1}, x^{i_0}y^{a_3}, \dots, x^{i_0}y^{a_{m-1}}), \\
A_2 &= (x^{i_0}y^{b_0}, x^{i_0}y^{b_2}, \dots, x^{i_0}y^{b_{m-2}}) \cup (x^{i_0}y^{b_1}, x^{i_0}y^{b_3}, \dots, x^{i_0}y^{b_{m-1}}),
\end{align*}
where each partition has at least $2$ distinct elements. This is possible since $S$ has at most $m-1$ identical elements. From Corollary \ref{q=1,r=2}, $A_1$ and $A_2$ each generate at least $q-1$ elements in $H$, therefore the product-set $A_1 \cdot A_2$ has to be $H$ by the Cauchy-Davenport inequality, which implies that $1 \in A_1 \cdot A_2$. Thus $S$ is not free of product-$1$ subsequences.

\item {\bf Subcase $q \equiv 3 \pmod 4$ and $q \ge 7$:} It is known that $-1$ is not a quadratic residue modulo $q$.  Let 
\begin{align*}
A_1 &= \{x^{i_0}y^{a_0}, x^{i_0}y^{a_1}, \dots, x^{i_0}y^{a_{m-1}}\}, \\
A_2 &= \{x^{i_0}y^{b_0}, x^{i_0}y^{b_1}, \dots, x^{i_0}y^{b_{m-1}}\},
\end{align*}
where $A_1$ and $A_2$ are disjoint subsequences of $S$ and each one has at least $2$ distinct elements. It is enough to consider the exponents of $y$ and construct the following sets of exponents:
\begin{align*}
X &= \{ s^{ni_0}(a_0 + a_1s^{i_0} + \dots + a_{m-2}s^{i_0(m-2)} + a_{m-1}s^{i_0(m-1)}) \in \Z_q | \; 0 \le n \le m-1 \}, \\
Y &= \{ s^{ni_0}(b_0 + b_1s^{i_0} + \dots + b_{m-2}s^{i_0(m-2)} + b_{m-1}s^{i_0(m-1)}) \in \Z_q | \; 0 \le n \le m-1 \}.
\end{align*}
Suppose that $0 \not\in X$, $0 \not\in Y$ and $0 \not\in X+Y$ (otherwise we are done). By Lemma \ref{pinA}, $|X| = |Y| = m$. Let 
$$\alpha = a_0 + a_1s^{i_0} + \dots + a_{m-2}s^{i_0(m-2)} + a_{m-1}s^{i_0(m-1)} \in X.$$
 If $-\alpha \in Y$ then $0 \in X+Y$, a contradiction, therefore $-\alpha \in X$. Hence, there exist $0 \le n \le m-1$ such that $\alpha s^{2ni_0} \equiv -\alpha \pmod q$, which implies $s^{2ni_0} \equiv -1 \pmod q$, so $-1$ is a quadratic residue modulo $q$, which is also a contradiction.
\end{enumerate}

\item {\bf Case $r = 1$:} In this case, $m = q-1$ and, in particular, $m$ is even and $s^{i_0}$ generates $\Z_q^*$. The cases where $(m,q) = (4,5)$ follow from Proposition \ref{mq45}, therefore we may assume $q \ge 7$. If there exist $m$ identical elements then their product is $1$, thus there are at most $m-1$ identical elements. Since 
$$\frac{m+q-2}{m-1} = \frac{2m-1}{m-1} > 2,$$
 there are at least $3$ distinct elements among $S = \{x^{i_0}y^{\alpha_1}, \dots, x^{i_0}y^{\alpha_{2m-1}}\}$. Again, we  split up this case into $q \equiv 1 \pmod 4$ and $q \equiv 3 \pmod 4$, as follows:

\begin{enumerate}[({3.}1)]
\item {\bf Subcase $q \equiv 1 \pmod 4$ and $q \ge 13$:} We may choose
$$A_1 = \{x^{i_0}y^{a_0}, x^{i_0}y^{a_1}, \dots, x^{i_0}y^{a_{m-1}}\}$$
 to be a subsequence of $S$ that we can further  split  into two subsequences of the same size $(q-1)/2$, say 
$$A_1 = \{x^{i_0}y^{a_0}, x^{i_0}y^{a_2}, \dots, x^{i_0}y^{a_{m-2}}\} \cup \{x^{i_0}y^{a_1}, x^{i_0}y^{a_3}, \dots, x^{i_0}y^{a_{m-1}}\},$$
 satisfying $a_0 \not\equiv a_2 \pmod q$ and $a_1 \not\equiv a_3 \pmod q$. Since $s^{i_0}$ generates $\Z_q^*$, $s^{2i_0}$ generates the quadratic residues of $\Z_q^*$. From Corollary \ref{q=1,r=1}, there exist two permutations $\sigma$ of $(0,2,\dots,m-2)$ and $\tau$ of $(1,3,\dots,m-1)$ such that
\begin{align*}
1 &\equiv a_{\sigma(0)} + a_{\sigma(2)}s^{2i_0} + \dots + a_{\sigma(2)}s^{(m-2)i_0} \pmod q, \\
-1 &\equiv a_{\tau(1)}s^{i_0} + a_{\tau(3)}s^{3i_0} + \dots + a_{\tau(m-1)}s^{(m-1)i_0} \pmod q.
\end{align*}
Therefore
\begin{align*}
1 = y \cdot y^{-1} &= y^{a_{\sigma(0)} + a_{\sigma(2)}s^{2i_0} + \dots + a_{\sigma(2)}s^{(m-2)i_0}}\cdot y^{a_{\tau(1)}s^{i_0} + a_{\tau(3)}s^{3i_0} + \dots + a_{\tau(m-1)}s^{(m-1)i_0}} \\
&= x^{i_0}y^{\sigma(0)} \cdot x^{i_0}y^{\tau(1)} \cdot x^{i_0}y^{\sigma(2)} \cdot x^{i_0}y^{\tau(3)} \dots x^{i_0}y^{\sigma(m-2)} \cdot x^{i_0}y^{\tau(m-1)},
\end{align*}
thus $S$ is not free of product-$1$ subsequences.

\item {\bf Subcase $q \equiv 3 \pmod 4$ and $q \ge 7$:} It is known that $-1$ is not a quadratic residue modulo $q$. Let 
$$A_1 = \{x^{i_0}y^{a_0}, x^{i_0}y^{a_1}, \dots, x^{i_0}y^{a_{m-1}}\},$$
 where at least $3$ of the $a_j$'s are distinct. By considering the set of all products obtained by changing the order of the elements of $A_1$, we obtain the set 
\begin{align*}
\mathcal C = \{\alpha_{\pi} = &\; a_{\pi(0)} + a_{\pi(1)}s^{i_0} + \dots + a_{\pi(m-1)}s^{(m-1)i_0} \in \Z_q | \\
&\; \pi \text{ is a permutation of }(0,1,\dots,m-1)\}.
\end{align*}
 Hence, it is enough to consider the exponents $\alpha_{\pi}$ of $y$. Reindexing the $a_j$'s, we may construct the set of exponents
\begin{align*}
X &= \{s^{2ni_0}(a_0 + a_2s^{2i_0} + \dots + a_{m-2}s^{i_0(m-2)}) \in \Z_q | \; 0 \le n \le m-1\}, \\
Y &= \{s^{2ni_0}(a_1s^{i_0} + a_3s^{3i_0} + \dots + a_{m-1}s^{i_0(m-1)}) \in \Z_q | \; 0 \le n \le m-1\},
\end{align*}
where $X$ and $Y$ each have  at least $2$ distinct elements. Clearly, $X+Y \subset \mathcal C$. Suppose that $0 \not\in X+Y$ (otherwise we are done). By Lemma \ref{pinA}, $|X+Y| = m$, thus $X+Y = \Z_q^*$. By Lemma \ref{wlog}, we may assume without loss of generality that $|X| = m/2$ and $|Y| = m/2$. Let 
$$\alpha = a_0 + a_2s^{2i_0} + \dots + a_{m-2}s^{i_0(m-2)} \in A.$$
 If $-\alpha \in Y$ then $0 \in X+Y$, a contradiction. Therefore $-\alpha \in X$. Hence, there exists $0 \le n \le m-1$ such that $\alpha s^{2ni_0} \equiv -\alpha \pmod q$, which implies $s^{2ni_0} \equiv -1 \pmod q$, so $-1$ is a quadratic residue modulo $q$, but this is impossible.
\end{enumerate}
\end{enumerate}

\qed

\Ack We would like to thank the anonymous ``user44191'' from the {\em Math Overflow Forum} (\url{http://mathoverflow.net/}) for providing a hint for solving Lemma \ref{sinvariant}. The second author would like to thank {\em CAPES}/Brazil for the PhD student fellowship.

\end{document}